\documentclass{amsart}
\usepackage{amsmath,amssymb,amsthm,amscd}
\theoremstyle{definition}

\theoremstyle{remark}

\begin{document}

\title
[Multivalued functionals, one-forms and deformed de Rham complex]
{Multivalued functionals, one-forms and  deformed de Rham complex}
\author{Dmitri V. Millionschikov}
\thanks{Partially supported by
the Russian Foundation for Fundamental Research, grant no.
05-01-01032} 
\address{Department of Mathematics and Mechanics, Moscow
State University, 119899 Moscow, RUSSIA}
\email{million@mech.math.msu.su}

\begin{abstract}
We discuss some applications of the Morse-Novikov theory to
some problems in modern physics, where appears a non-exact closed
$1$-form $\omega$(multi-valued functional).
We focus mainly our attention to
the cohomology $H^*_{\lambda \omega}(M^n, {\mathbb R})$ of
the de Rham complex $\Lambda^*(M^n)$ of a compact manifold $M^n$
with a deformed differential $d_{\omega}=d +\lambda \omega$.
Using Witten's approach to the Morse theory one can estimate the number of
critical points of $\omega$ in terms of
$H^*_{\lambda \omega}(M^n, {\mathbb R})$ with sufficiently large values of
$\lambda$ (torsion-free Novikov's inequalities).

We show that for an interesting class of solvmanifolds
$G/\Gamma$ the cohomology $H^*_{\lambda \omega}(G/\Gamma, {\mathbb R})$
can be computed as the cohomology $H^*_{\lambda \omega}(\mathfrak{g})$ of
the corresponding Lie algebra $\mathfrak{g}$ associated with the
one-dimensional representation $\rho_{\lambda \omega}$.
Moreover $H^*_{\lambda \omega}(G/\Gamma, {\mathbb R})$
is almost always trivial except a finite number of classes
$[\lambda \omega]$ in $H^1(G/\Gamma, {\mathbb R})$.
\end{abstract}

\date{}

\maketitle

\section*{Introduction}
In the begining of the 80-th S.P. Novikov  constructed
(\cite{N1}, \cite{N2})
an analogue of the Morse
theory for smooth multi-valued functions, i.e. smooth
closed $1$-forms
on a compact smooth manifold $M$.
In particular he introduced the Morse-type inequalities
(Novikov's inequalities)
for numbers $m_p(\omega)$ of zeros
of index $p$ of a Morse $1$-form $\omega$.

In ~\cite{N3}, ~\cite{Pa}  a method
of obtaining the torsion-free Novikov
inequalities in terms of the de Rham complex of manifold was proposed.
This method was based on Witten's approach ~\cite{W} to the Morse theory.
A. Pazhitnov generalized Witten's deformation $d+tdf$ ($f$ is a Morse
function on $M$) of the standard
differential $d$ in $\Lambda^*(M)$ by replacing $df$ by
an arbitrary Morse $1$-form on $M$.
For sufficiently large real values $t$
one have the following estimate (torsion-free Novikov
inequality ~\cite{Pa}):
$$ m_p(\omega) \ge \dim H^p_{t \omega}(M,{\mathbb R}),$$
where  by  $H^p_{t \omega}(M,{\mathbb R})$ we denote the $p$-th cohomology
group of the de Rham complex $(\Lambda^*(M), d+t \omega)$
with respect to the new deformed differential
$d+t \omega$.

Taking a  complex parameter $\lambda$
one can identify  $H^*_{\lambda \omega}(M,{\mathbb C})$
with the cohomology $H^*_{\rho_{\lambda \omega}}(M_n,{\mathbb C})$
with coefficients in the local system
$\rho_{\lambda \omega}$ of groups ${\mathbb C}$, where
$$\rho_{\lambda \omega}(\gamma)=\exp{\int_{\gamma} \lambda \omega},
\gamma \in \pi_1(M).$$

L. Alania in ~\cite{Al} studied
$H^*_{\rho_{\lambda \omega}}(M_n,{\mathbb C})$ of a
class of nilmanifolds $M_n$. He proved that
$H^*_{\rho_{\lambda \omega}}(M_n,{\mathbb C})$ is trivial if
$\lambda \omega \ne 0$.
The proof was based on the Nomizu theorem ~\cite{Nz} that reduce the problem
to the computations in terms of the corresponding nilpotent Lie algebra.
It was remarked in \cite{Mil} that triviality of
$H^*_{\rho_{\lambda \omega}}(G/\Gamma,{\mathbb R})$, with
$\lambda \omega \ne 0$
follows from Dixmier's theorem ~\cite{D}, namely:
{\it for a nilmanifold} $G/\Gamma$ {\it the cohomology
$H^*_{\omega}(G/\Gamma, {\mathbb R})$ coincides with
the cohomology $H^*_{\omega}(\mathfrak{g})$
associated with the one-dimensional representation of the Lie algebra}
$\rho_{\omega}: \mathfrak{g} \to {\mathbb R},
\rho_{\omega }(\xi) = \omega(\xi)$ {\it and
hence} $H^*_{\omega}(G/\Gamma, {\mathbb R})=
H^*_{\omega}(\mathfrak{g})=0$.

Applying Hattori's theorem ~\cite{H}
one can observe that the isomorphism
$$H^*_{\omega}(G/\Gamma, {\mathbb R})
\cong H^*_{\omega}(\mathfrak{g})$$ still holds on for
compact solvmanifolds $G/\Gamma$ with completely solvable Lie
group $G$. The calculations
show that the cohomology $H^*_{\omega}(G/\Gamma, {\mathbb R})$
can be non-trivial for certain values
$[\omega] \in H^1(G/\Gamma, {\mathbb R})$. However there exist
only a finite number of such values.

Let us consider a finite subset
$\Omega_{G/\Gamma}$ in $H^1(G/\Gamma, {\mathbb R}) \cong
H^1(\mathfrak{g})$:
$$
\Omega_{G/\Gamma}=\left\{ \alpha_{i_1}{+} \dots{+} \alpha_{i_s} |
\;\; 1 \le i_1{<} \dots {<} i_s \le n, \;\; s=1, {\dots}, n \right\},
$$
where the set $\{\alpha_{1}, \dots, \alpha_{n} \}$ of closed
$1$-forms is in fact the set of the weights of
completely reducible representation associated to the adjoint
representation of $\mathfrak{g}$.
It was proved in \cite{Mil}:
{\it if} $-[\omega] \notin \Omega_{G/\Gamma}$,
{\it then the cohomology} $H^*_{\omega}(G/\Gamma, {\mathbb R})$
{\it is trivial.}

\section{Dirac monopole, multivalued actions and Feinman quantum amplitude}
The notion of multivalued functional originates from topological
study of the quantization process of the motion of a charged particle
in the field of a Dirac monopole \cite{Dir}. The Kirchhoff-Thomson
equations for free motion of solids in a perfect
noncompressible liquid also can be reduced
to the theory of a charged particle on the sphere $S^2$ with some metric
$g_{\alpha\beta}$
in a potential field ${\it U}$ and in an effective magnetic field $F=F_{12}$
with a non-zero flux $4\pi s$ through $S^2$. Locally (in some domain
$U_{\alpha}$) on the sphere we have the following formula
for the action $S_{\alpha}(\gamma)$:
\begin{equation}
S_{\alpha}(\gamma)=\int_{\gamma}\left(\frac{1}{2}g_{ij}
{\dot x}^{i}{\dot x}^{j}-{\it U}+e{\it A}^{\alpha}_k{\dot x}^{k}
\right)dt,
\end{equation}
where
\begin{equation}
x^1=\theta, x^2=\varphi,
\quad F_{12}d\theta {\wedge} d\varphi= d({\it A}^{\alpha}_k{d x}^{k}),
\quad \int \int_{S^2}F_{12}d\theta {\wedge} d\varphi =4 \pi s \ne 0.
\end{equation}

One can consider Feynman's paths integral approach to the
quantization of the problem considered above. Recall that
in the standard situation of single-valued action $S$, we consider
the amplitude $$\exp{\{2\pi iS(\gamma)\}},
\quad \gamma \in \Omega(x,x')$$
and the propagator
$$K(x,x')=\int_{\Omega(x,x')}\exp{\{2\pi iS(\gamma)\}}
D\gamma.$$
For the Dirac monopole one can consider the set $\{S_1, S_2\}$
of local actions where
$U_1=S^2\backslash P_N$ and  $U_2=S^2\backslash P_S$, by
$P_N, P_S$ we denote the poles of the sphere $S^2$.
Taking the equator $\gamma$ with the positive  orientation, one
can easily test the ambiguity of the action:
\begin{equation}
S_1(\gamma)-S_2(\gamma)=
e \int_{\gamma}(A^1_kdx^k-A^2_kdx^k)=
e \int \int_{S^2}F_{12}d\theta {\wedge} d\varphi =4 \pi s e \ne 0.
\end{equation}
The monopole is quantized if and only if the amplitude
$\exp{\{2\pi iS_{\alpha}(\gamma)\}}$
is a single-valued functional, i.e. for an arbitrary
closed $\gamma \in U_1 \cap U_2$ we have
$$\exp{\{2\pi iS_1(\gamma)\}}=
\exp{\{2\pi iS_2(\gamma)\}}$$
The last condition is equivalent
to the following one:
\begin{equation}
4 \pi s e = k, \; k \in {\mathbb Z}.
\end{equation}

Generalizing the situation with the Dirac monopole Novikov \cite{N2}
considered a $n$-dimensional manifold $M^n, n >1$ with a metric
$g_{ij}$, with a scalar potential $U$ and with a two-form $F$ of magnetic
field not necessarily exact.
In these settings one can consider a set of open $U_{\alpha} \subset M^n$,
such that  $F=F_{ij}dx^i {\wedge} dx^j$ is exact on $U_{\alpha}$ and
$M^n \subset \cup_{\alpha} U_{\alpha}$. A $1$-form
$\omega_{\alpha}= A_{k}^{\alpha}{d x}^{k}$,
$d\omega_{\alpha}=F_{ij}dx^i {\wedge} dx^j$
is determined up to some closed $1$-form and we can consider
the set of local actions:
\begin{equation}
S_{\alpha}(\gamma)=\int_{\gamma}\left(\frac{1}{2}g_{ij}
{\dot x}^{i}{\dot x}^{j}-{\it U}+e{\it A}^{\alpha}_k{\dot x}^{k}
\right)dt,
\end{equation}

Let us consider a path $\gamma \subset U_{\alpha} \cap U_{\beta}$.
The values $S_{\alpha}(\gamma)$ and $S_{\beta}(\gamma)$
do not coincide generally speaking. Hence the set $\{S_{\alpha}\}$
of local actions defines a multi-valued functional $S$.
As $\omega_{\alpha}-\omega_{\beta}$ is closed on
$U_{\alpha} \cap U_{\beta}$ the integral
$$S_{\alpha}(\gamma_{\lambda})-S_{\beta}(\gamma_{\lambda})=
\int_{\gamma_{\lambda}}(\omega_{\alpha}-\omega_{\beta})$$
is invariant under any deformation
$\gamma_{\lambda} \subset U_{\alpha} \cap U_{\beta}$ of
$\gamma$ in the class:
a) periodic curves;
b) the curves with the same end-points.

The crucial Novikov's observation was:
the infinite-dimensional $1$-form $\delta S$ is well-defined
and closed for the following function spaces:

a) $\Omega^+$ of the oriented
closed contours $\gamma$, such that $\exists \alpha, \gamma \subset
U_{\alpha}$;

b) $\Omega(x,x')$ of the paths $\gamma(x,x')$ joining points $x, x'$,
such that $\exists \alpha, \gamma(x,x') \subset U_{\alpha}$.

The set $\{ S_{\alpha}\}$
of local actions determines also a multivalued (in general)
functional $\exp{\{2\pi iS\}}$
on $\Omega(x,x')$. The local variational system  $\{ S_{\alpha}\}$
is quantized if and only if the Feinman quantum amplitude
$\exp{\{2\pi iS\}}$ is a single-valued functional on $\Omega^+$.
Or, in other words, for all $\gamma \in U_{\alpha} \cap U_{\beta}$ we have
\begin{equation}
\label{int_int}
\int_{\gamma}(\omega_{\alpha}-\omega_{\beta})= k,\quad k \in {\mathbb Z}
\end{equation}
If $U_{\alpha}$ and $U_{\beta}$ are simply connected domains in $M^n$
it is possible to consider a map $f: S^2 \to M$ such that $\gamma$
is the image
of the equator of the sphere $S^2$ and the images of two half-spheres of
$S^2$ lie in $U_{\alpha}$ and $U_{\beta}$ respectively. Then the condition
(\ref{int_int}) can be rewrited as
\begin{equation}
\int_{f(S^2)}F_{ij}dx^i \wedge dx^j= k,\quad k \in {\mathbb Z}.
\end{equation}
Hence we can propose the following sufficient
condition of the quantization:
{\it a local variational system $\{S_{\alpha}\}$ on $M^n$ that corresponds
to some magnetic field $F=F_{ij}dx^i \wedge dx^j$ is quantized if
$F$ has integer-valued fluxes through all basic cycles
of $H_2(M^n,{\mathbb Z})$.}

One can remark that the last condition is in fact excessive:
it is sufficient to require integer-valued integrals of $F$
over spheric cycles that lie in the image of the Hurevich map
$H: \pi_2(M^n) \to H_2(M^n, {\mathbb Z})$.

\section{Aharonov-Bohm field and equivalent quantum systems}
Another interesting example
comes from Aharonov-Bohm experiment. We consider
the electron move outside the ideal endless solenoid, i.e. the
configuration space is $M=({\mathbb R}^2\backslash D_{\varepsilon})\times
{\mathbb R}=\{(x,y,z) \in {\mathbb R}^3, x^2{+}y^2 {>} \varepsilon^2\}$,
where $D_{\varepsilon}$ is two-dimensional disk of radius
$\varepsilon \to 0$. The magnetic field $F=F_{ij}dx^i \wedge dx^j$
vanish outside solenoid,
i.e. $F\equiv 0$ on $M$, hence
\begin{equation}
S_{\omega_{\alpha}}(\gamma)=\int_{\gamma}\frac{m{\dot x}^2}{2}dt
+\omega_{\alpha},
\end{equation}
where $\omega_{\alpha}=e A_kdx^k$ is an arbitrary
closed $1$-form on $M$.
The cohomology space
$$H^1(M, {\mathbb R})=
H^1({\mathbb R}^2\backslash D_{\varepsilon}, {\mathbb R})=
H^1(S^1, {\mathbb R})={\mathbb R}$$ is one-dimensional
and hence
$$\omega_{\alpha}=\frac{e \Phi_{\alpha}}{2\pi }
\frac{xdy-ydx}{x^2+y^2}+df_{\alpha},$$
for some constant $\Phi_{\alpha}$ and function $f_{\alpha}$ on $M$.

Taking the circle $\gamma_0=\partial D_{\varepsilon}=
\{(\varepsilon\cos{\varphi}, \varepsilon\sin{\varphi},0),
0\le \varphi < 2\pi\}$
we have
$$\int_{\gamma_0} A_k^{\alpha}dx^k=
\frac{1}{e}\int_{\gamma_0} \omega_{\alpha}=\Phi_{\alpha}=
\int_{D_{\varepsilon}}F_{12}dx \wedge dy.$$
Hence the constant $\Phi_{\alpha}$ is equal to the flux of the magnetic field
$F$ through the orthogonal section $D_{\varepsilon}$ of our solenoid.

The form $\omega_{\alpha}$ determines a representation
$\rho_{\omega_{\alpha}}$ of the fundamental
group of $M$:
$$\rho_{\omega_{\alpha}}: \pi_1(M) \to {\mathbb C}^*, \quad
\rho_{\omega_{\alpha}}(\gamma)=
\exp{\{2\pi i\int_{\gamma}\omega_{\alpha}\}}, \;\;
\gamma \in \pi_1(M).$$

Let $S_{\omega_1}$ and $S_{\omega_2}$ be two actions for our system.
They are quantummechanically equivalent if and only if
$$\exp{\{2\pi iS_{\omega_1}(\gamma)\}}=c(x,x')
\exp{\{2\pi iS_{\omega_2}(\gamma)\}},$$
with a phase factor $c(x,x')$ depending only on end points $x, x'$
of $\gamma$ and $|c(x,x')|=1$, i.e. $c(x,x')$ is physically unobservable.
It is easy to show that the actions
$S_{\omega_1}$ and $S_{\omega_2}$
are quantummechanically equivalent if and only if
for any loop $\gamma \in \pi_1(M)$ the value of the integral
$\int_{\gamma}(\omega_1-\omega_2)$ is integer or, in other
words, the form $(\omega_1-\omega_2)$ has
integer-valued integrals over basic cycles of $H_1(M, {\mathbb Z})$.

In our case
$H_1(M, {\mathbb Z})={\mathbb Z}$ and the last condition is equivalent to
the following one
\begin{equation}
\int_{\gamma_0}(\omega_1-\omega_2)=
e(\Phi_1-\Phi_2)=k, \quad k \in {\mathbb Z}.
\end{equation}
Hence (one of the important observations in the Aharonov-Bohm experiment)
the fields with fluxes $\Phi_1$ and $\Phi_2$, such that
$\Phi_1-\Phi_2=\frac{k}{e}, k \in {\mathbb Z}$ can not be
distinguished by any interference effect.

Now let us consider the case when $M^n$ is not simply connected
and the two-form $F$ is globally exact on $M^n$ (like in the
Aharonov-Bohm experiment). Two sulutions $\omega_1, \omega_2$
of the equation $d\omega=F_{ij}dx^i \wedge dx^j$ that correspond
to two different actions $S_1(\gamma)$ and $S_2(\gamma)$ are
determined up to a differential $df$ by their integrals
$\int_{\gamma_k}\omega_i$ over
the basic cycles $\gamma_k$ of $H_1(M^n, {\mathbb Z})$. These integrals
can be interpreted as the fluxes of the continuation of $F$
(with possible singularities) to some large manifold $\tilde M^n$.
Two variational systems $S_1(\gamma)$ and $S_2(\gamma)$ are
quantummecanically equivalent if and only if all
integrals $\int_{\gamma_k}(\omega_1-\omega_2)$  over basic
cycles $\gamma_k$ of $H_1(M^n, {\mathbb Z})$ are integer-valued.

The form $\omega_{12}=\omega_1-\omega_2$ is a closed $1$-form
on $M^n$ and it determines
a representation $\rho_{\omega_{12}}$ of the fundamental group $\pi_1(M^n)$:
$$\rho_{\omega_{12}}: \pi_1(M^n) \to {\mathbb C}^*, \quad
\rho_{\omega_{12}}(\gamma)=
\exp{\{2\pi i\int_{\gamma}\omega_{12}\}}, \;\;
\gamma \in \pi_1(M^n).$$

Let $M$ be a finite-dimensional (or infinite-dimensional) manifold
and $S: M \to {\mathbb R}$ a function (functional) on it.

What are the relations between the set of the stationary points $dS=0$
($\delta S=0$) and the topology of the
manifold $M$ ?

If $S$ is a Morse function (generic situation),
i.e. $d^2S$ is non-degenerate at critical points,
then one can define the Morse index $ind(P)$ of a critical point
$P$ of $S$ as the number of negative squares of the quadratic form
$d^2S(P)$ (if it is finite in the infinite-dimensional case).

Under some natural hypotheses the following inequality
can be established:
$$m_p(S) \ge b_p(M)=\dim{H^p(M)}.$$

\begin{figure}
\begin{picture}(80,30)(-15,2)
  \put(30,0){\line(0,1){30}}
  \put(30,20){\vector(0,1){10}}
  \put(10,15){\circle{7}}
  \put(10,15){\circle{20}}
  \put(10,8){\circle*{1}}
  \put(10,11.5){\circle*{1}}
  \put(10,18.5){\circle*{1}}
  \put(10,22){\circle*{1}}
  \put(30,8){\circle*{1}}
  \put(30,11.5){\circle*{1}}
  \put(30,18.5){\circle*{1}}
  \put(30,22){\circle*{1}}
  \multiput(12,8)(2,0){9}{\line(1,0){1}}
  \multiput(12,11.5)(2,0){9}{\line(1,0){1}}
  \multiput(12,18.5)(2,0){9}{\line(1,0){1}}
  \multiput(12,22)(2,0){9}{\line(1,0){1}}
  \put(33,7){$ind(q_1){=}0$}
  \put(33,11.5){$ind(q_2){=}1$}
  \put(33,18){$ind(q_3){=}1$}
  \put(33,22){$ind(q_4){=}2$}
  \put(17.5,28.5){$f(q){=}z$}
  \put(26.5,6){$z_1$}
  \put(26.5,13){$z_2$}
  \put(26.5,19.5){$z_3$}
  \put(26.5,23){$z_4$}
  \put(8,5){$q_1$}
  \put(9.5,13.5){$q_2$}
  \put(7.7,19.8){$q_3$}
  \put(10,24){$q_4$}
  {\linethickness{0.5mm}
  \put(53,0){\line(0,1){30}}}
  \put(55,23){$m_2(f)=\dim{H^2({\mathbb T}^2,{\mathbb R})=1}$}
  \put(55,15){$m_1(f)=\dim{H^1({\mathbb T}^2,{\mathbb R})=2}$}
  \put(55,7){$m_0(f)=\dim{H^0({\mathbb T}^2,{\mathbb R})=1}$}
                 \end{picture}
\caption{A height-function $f(q){=}z$ for ${\mathbb T}^2$.}
\label{fig0}
\end{figure}
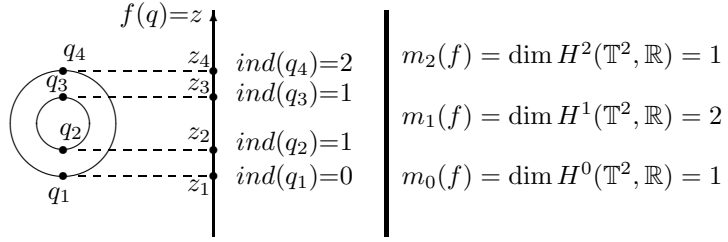

\section{Semiclassical motion of electron and
critical points of $1$-form}

The semiclassical model of electron motion is an important tool
for investigating conductivity in crystals under the action of a magnetic
field. In the same time it is one of the most important
examples of applications of topological methods in the modern physics.

Let us consider the corresponding quantum
system defined for some crystal lattice $L={\mathbb Z}^3$.
Its eigenstates are the Bloch functions $\psi_p$. The particle
quasimomentum $p$ is defined up to a vector of the dual lattice
$L^*={\mathbb Z}^3$. Hence one can regard the space of quasimomenta
as $3$-dimensional torus
${\mathbb T}^3={\mathbb R}^3/{\mathbb Z}^3$.
The state energy $\varepsilon(p)$ is thus a function on ${\mathbb T}^3$,
i.e. $3$-periodical function in ${\mathbb R}^3$.

An external homogeneous constant magnetic field
is a constant vector $H=(H_1,H_2,H_3)$ or in other words it is
a $1$-form $\omega=H_1dp_1+H_2dp_2+H_3dp_3$ with constant coefficients.

The semiclassical trajectories projected to the space of quasimomenta
are connected components of the intersection of the planes
$(p,H)=const$ with constant energy surfaces $\varepsilon(p)=const$.

The constant energy surfaces $\varepsilon(p)=\varepsilon_F$
that correspond to the Fermi energies $\varepsilon_F$ are called
the Fermi surfaces. There are nonclosed trajectories on the Fermi surfaces
with asymtotic directions and this topological fact
explains an essential anisotropy of the metal conductivity
at low temperatures.

One can study the intersections
$$(p,H)=c_0, \quad \varepsilon(p)=\varepsilon_0$$
as the level surfaces of the $1$-form
$$\omega_{\varepsilon_0}=
(H_1dp_1+H_2dp_2+H_3dp_3)|_{{\hat M}_{\varepsilon_0}},$$
where  $2$-dimensional manifold
$${\hat M}_{\varepsilon_0}=
\{p \in {\mathbb R}^3|\varepsilon(p)=\varepsilon_0\}$$
is the universal covering
of the compact Fermi surface $\varepsilon(p)=\varepsilon_0$
in ${\mathbb T}^3$. The last one we will denote also by $M_{\varepsilon_0}$.
We can treat the $3$-periodic form $\omega_{\varepsilon_0}$ as
a $1$-form on the compact manifold $M_{\varepsilon_0}$
(we will keep the same notation for it).

The information about critical points of $\omega_{\varepsilon_0}$
is very important in the problem considered above. A generic
$1$-form $\omega_{\varepsilon_0}$ is a Morse form and has finitely many
critical points on $M_{\varepsilon_0}$.

\section{Witten's deformation of de Rham complex and Morse-Novikov theory}
In 1982 E.~Witten
proposed a new beautiful proof of the Morse inequalities
using some analogies with supersymmetry quantum mechanics.
Taking an arbitrary smooth
real-valued function $f$ on a Riemannian manifold $M^n$
he considered a new deformed differential $d_t$
in the de Rham complex $\Lambda^*(M^n)$ ($t$ is a real parameter):
\begin{equation}
\begin{split}
d_t=e^{{-}ft}de^{ft}=d+tdf\wedge, \\
 d_t(\xi)=d\xi+tdf\wedge \xi, \; \xi \in \Lambda^*(M^n),
\end{split}
\end{equation}
where $d$ is the standard differential in $\Lambda^*(M^n)$:
\begin{equation}
\begin{split}
d: \Lambda^p(M^n)\to \Lambda^{p{+}1}(M^n),\\
\xi=\sum\limits_{i_1{<}\dots{<}i_p}\xi_{i_1{\dots}i_p}
dx^{i_1}{\wedge}\dots{\wedge}dx^{i_p} \; \in \Lambda^p(M^n),\\
d\xi=\sum\limits_{i_1{<}\dots{<}i_p}
\sum\limits_q
\frac{\partial \xi_{i_1{\dots}i_p}}{\partial x^q}
dx^q{\wedge}
dx^{i_1}{\wedge}\dots{\wedge}dx^{i_p}\; \in \Lambda^{p{+}1}(M^n).
\end{split}
\end{equation}
Taking arbitrary smooth vector fields $X_1,\dots, X_{p{+}1}$
on $M^n$ we have also the following formula:
\begin{equation}
\label{Cartan_form}
\begin{split}
d\xi(X_1, {\dots}, X_{p{+}1})=
\sum_{1{\le}i{<}j{\le}p{+}1}({-}1)^{i{+}j}
\xi([X_i,X_j],X_1, {\dots}, \hat X_i,
{\dots}, \hat X_j, {\dots}, X_{q{+}1})+\\
+ \sum_{i}({-}1)^{i{+}1}X_i \xi(X_1, {\dots}, \hat X_i,
{\dots}, X_{p{+}1}).
\end{split}
\end{equation}
We recall that that a differential $p$-form $\xi$ is called closed
if $d\xi=0$ and it is called exact if $\xi=d\xi'$ for some
$(p{-}1)$-form $\xi'$. As $d^2=0$ the space of exact forms is
a subspace of the space of closed ones and the $p$-th de Rham
cohomology group $H^p(M^n,{\mathbb R})$ of the manifold $M^n$
is defined as a quotient space of closed $p$-forms
modulo exact ones. In the same manner the cohomology
$H_t^*(M^n,{\mathbb R})$ of the de Rham complex with respect to
the deformed differential $d_t$ can be defined.

The operators $d_t$ and $d$ are conjugated by the invertible
operator $e^{ft}$ and therefore the cohomology groups
$H^*(M^n,{\mathbb R})$ (the standard de Rham cohomology) and
$H_t^*(M^n,{\mathbb R})$ (the new ones) are isomorphic to each other.
On the level of the forms this isomorphism is given by the gauge
transformation
$$\xi \to e^{ft}\xi.$$

One can define the adjoint operator
$d_t^*=e^{ft}d^*e^{{-}ft}$ with respect to the scalar product of differential forms
$$(\alpha,\beta)=\int_{M^n}(\alpha,\beta)_xdV,$$
where $(\alpha,\beta)_x$ is a scalar product in the bundle
$\Lambda^*(T_x^*(M^n))$ evaluated with respect to the
Riemannian metric $g_{ij}$ of $M^n$
and $dV$ is the corresponding volume form.

One can also consider the deformed Laplacian
$H_t=d_t d_t^* + d_t^*d_t$ acting on forms.
An arbitrary element $\omega$ from $H_t^p(M^n,{\mathbb R})$
can be uniquely represented
as an eigenvector
with zero eigenvalue of the Hamiltonian $H_t=d_t d_t^* + d_t^*d_t$.
Hence one can compute the Betti number $b_p(M^n)=\dim{H^p(M^n,{\mathbb R})}$
as the number of zero eigenvalues of $H_t$ acting on $p$-forms.

It can be calculated that
\begin{equation}
H_t=d_t d_t^* + d_t^*d_t=d d^* + d^*d+t^2(df)^2+
t\sum\limits_{i{,}j}\nabla^2_{(i{,}j)}(f)[\tilde a^i, \tilde a^{j*}],
\end{equation}
where $(df)^2=(df,df)_x=g^{ij}\frac{\partial f}{\partial x^i}
\frac{\partial f}{\partial x^j}$ and
$$\tilde a^i(\xi)=dx^i \wedge \xi, \quad \nabla^2_{(i{,}j)}=
\nabla_i \nabla_j -\Gamma^k_{ij}\nabla_k.$$

As the "potential energy" $t^2(df)^2$ of the Hamiltonian $H_t$
becomes very large for $t \to +\infty$ the eigenfuctions
of $H_t$ are concentrated near the critical points  $df=0$
and the low-lying eigenvalues of $H_t$ can be calculated by expanding
about the critical points. Taking the Morse coordinates
$x^i$ in some neighbourhood $W$ of a critical point $P$
$$f(x)=\frac{1}{2}\sum\lambda_i(x^i)^2, \quad
\lambda_1{=}\dots{=}\lambda_q{=}{-}1,
\lambda_{q{+}1}{=}\dots{=}\lambda_n{=}1,$$
where $q$ is the index of the critical point $P$,
and introducing a Riemmanian metric $g_{ij}$ on $M^n$ such that
$x^i$ are Euclidean coordinates for $g_{ij}$ in $W$
one can locally evaluate the Hamiltonian $H_t$:
\begin{equation}
H_t=\sum\limits_i\left({-}\frac{\partial^2}{\partial{x^i}^2}
+t^2{x^i}^2+t\lambda_i[\tilde a^i, \tilde a^{i*}]\right).
\end{equation}
The operator
$$H_i={-}\frac{\partial^2}{\partial{x^i}^2}+t^2{x^i}^2$$
is the Hamiltonian of the simple harmonic oscillator and it has the
following set of eigenvalues
$$t(1+2N_i), \quad N_i=0,1,2,\dots$$
with simple multiplicities. The operator
$H_i$ commutes with $[\tilde a^i, \tilde a^{i*}]$
and the eigenvalues of the last operator are equal to $\pm 1$:
$$
[\tilde a^i, \tilde a^{i*}]
(\psi (x)dx^{i_1}{\wedge}\dots{\wedge}dx^{i_p})=
\left\{ \begin{array}{c}
\psi(x) dx^{i_1}{\wedge}\dots{\wedge}dx^{i_p}, \quad i \in (i_1,\dots,i_p), \\
-\psi(x) dx^{i_1}{\wedge}\dots{\wedge}dx^{i_p}, \quad i \notin
(i_1,\dots,i_p).
\end{array}\right. $$
Hence the eigenvalues of the restriction $H_t|_W$
are equal to
\begin{equation}
t\sum\limits_i(1+2N_i+\lambda_i l_i),
\quad N_i=0,1,2,\dots,\:
l_i=\pm 1.
\end{equation}

The corresponding eigenfunctions
$\Psi_t=\psi(x,t) dx^{i_1}{\wedge}\dots{\wedge}dx^{i_p}$ are defined
in $W$ and not on the whole manifold $M^n$.
Using the partition of unit one can define a new smooth $q$-form
$\tilde \Psi_t$ on $M^n$ such that $\tilde \Psi_t$ coincides with
$\Psi_t$ in some $\tilde W \subset W$ and $\tilde \Psi_t \equiv 0$ outside
of $W$. The $q$-form $\tilde \Psi_t$ is called a quasi-mode:
\begin{equation}
\label{asympt_v}
H_t\tilde \Psi_t=t\left(\sum\limits_i(1+2N_i+\lambda_i l_i)+\frac{B}{t}
+\frac{C}{t^2}+\dots \right) \tilde \Psi_t, \quad t\to +\infty.
\end{equation}
The numbers $t\sum\limits_i(1+2N_i+\lambda_i l_i)$ are called
asymtotic eigenvalues and the minimal value $E_0^{as}$ of them approximates
the minimal eigenvalue of $H_t$ as $t \to +\infty$.

In order to find $E_0^{as}$, we must set
$N_i=0$
for all $i$. The sum
$$
\sum\limits_{i=1}^q(1 - l_i)+\sum\limits_{i{=}q{+}1}^n(1 + l_i).
$$
is non-negative and it is equal to zero if and only if
$$l_1=\dots=l_q=1, \quad l_{q{+}1}=\dots=l_n={-}1.$$
This means that, $H_t$ has precisely one zero asymtotic eigenvalue
for each critical point of index $q$. Hence we have precisely
$m_q(f)$ asymptotic zero eigenvalues (for $q$-forms).
Vanishing of the first term
of the asymtotical expansion (\ref{asympt_v}) for a minimal eigenvalue
of $H_t$ is only a necessary condition to have zero energy level,
hence the number $b_q(M^n)$ of zero eigenvalues does not exceed
the number of zero asymtotic eigenvalues.
In other words we have established the Morse inequalities
$$m_q(f) \ge b_q(M^n).$$

It was A.~Pajitnov who remarked that it is possible to apply
Witten's approach to the Morse-Novikov theory ~\cite{Pa}.
Let $\omega$ be a closed $1$-form on $M^n$
and $t$ a real parameter. As in the construction above
one can define a new
deformed differential $d_{t \omega}$ in $\Lambda^*(M)$
$$d_{t \omega}=d+t \omega\wedge, \quad
d_{t \omega}(\xi)=da+t \omega \wedge \xi.$$
If the $1$-form $\omega$ is not exact,
the cohomology $H^*_{t \omega}(M,{\mathbb R})$ of the de Rham
complex with the deformed differential $d_{t \omega}$
generally speaking is not isomorphic to the standard one $H^*(M,{\mathbb R})$.
But $H^*_{t \omega}(M,{\mathbb R})$ depends only on the
cohomology class of $\omega$:
for any pair $\omega, \omega'$ of $1$-forms such that
$\omega- \omega'=d \phi$, where $\phi$ is a smooth function on $M^n$
the cohomology
$H^*_{t \omega}(M,{\mathbb R})$ and
$H^*_{t \omega'}(M,{\mathbb R})$
are isomorphic to each other. This
isomorphism can be given by the gauge transformation
$$
\xi \to e^{t \phi}\xi; \quad
d \to e^{t \phi} d e^{-t \phi}=
d+ t d \phi \wedge.
$$
It is convinient also to consider instead of $t$ a complex
parameter $\lambda$.
It was remarked in ~\cite{N3}, ~\cite{Pa} that the cohomology
$H^*_{\lambda \omega}(M,{\mathbb C})$
of $\Lambda^*(M)$ with respect to the deformed differential
$d_{\lambda \omega}$ coincides with the cohomology
$H^*_{\rho_{\lambda \omega}}(M,{\mathbb C})$
with coefficients in the representation
$\rho_{\lambda \omega}: \pi_1(M) \to {\mathbb C}^*$ of fundamental
group defined by the formula
$$
\rho_{\lambda \omega}([\gamma])=\exp{\int_{\gamma} \lambda \omega},
\quad [\gamma] \in \pi_1(M),
$$

We denote corresponding Betti numbers by
$b_p(\lambda, \omega)$,
$b_p(\lambda, \omega)=\dim H^*_{\rho_{\lambda \omega}}(M,{\mathbb C})$.

There is another interpretation of
$H^*_{\rho_{\lambda \omega}}(M,{\mathbb C}$:
the representation $\rho_{\lambda \omega}: \pi_1(M) \to
{\mathbb C}^*$
defines a local system of groups ${\mathbb C}^*$
on the manifold $M$.
The cohomology of $M$
with coefficients in this local system coincides
with $H^*_{\rho_{\lambda \omega}}(M, {\mathbb C})$.

Now we can assume that $\omega$ is a Morse $1$-form, i.e., in a neighborhood
of any point $\omega=df$, where $f$ is a Morse function.
In other words $\omega$ gives a multi-valued Morse function.
The zeros of $\omega$ are isolated, and one can
define the index of each zero. The number of zeros of $\omega$
of index $p$ is denoted by $m_p(\omega)$.

Following Witten's scheme A. Pazhitnov showed in \cite{Pa} that
for sufficiently large real numbers $\lambda$
$$ m_p(\omega) \ge b_p(\lambda, \omega).$$

\section{Solvmanifolds and left-invariant forms}

A solvmanifold (nilmanifold) $M$
is a compact homogeneous space
of the form
$G/\Gamma,$ where $G$ is a simply connected solvable (nilpotent) Lie group
and $\Gamma$ is a lattice in $G$.

Let us consider some examples of solvmanifolds
(first two of them are nilmanifolds):

1) a $n$-dimensional torus $T^n={\mathbb R}^n/{\mathbb Z}^n$;

2) the Heisenberg manifold $M_3={\mathcal H}_3/\Gamma_3$, where
${\mathcal H}_3$ is the group of all matrices of the form
$$
   \left( \begin{array}{lcr}
   1 & x & z\\
   0 & 1 & y\\
   0 & 0 & 1\\
   \end{array} \right) , ~~~ x,y,z \in {\mathbb R},
$$
and a lattice $\Gamma_3$ is a subgroup of matrices with integer entries
$x,y,z \in {\mathbb Z}$.
$$
e_1=\begin{pmatrix}
0&1&0\\
0&0&0\\
0&0&0
\end{pmatrix},\quad e_2=
\begin{pmatrix}
0&0&0\\
0&0&1\\
0&0&0
\end{pmatrix}, \quad e_3=
\begin{pmatrix}
0&0&1\\
0&0&0\\
0&0&0
\end{pmatrix},
$$
and the only one non-trivial structure relation:
$[e_1,e_2]=e_3$.
The left invariant $1$-forms on ${\mathcal H}_3$
\begin{equation}
e^1=dx, \quad e^2=dy, \quad e^3=dz -xdy
\end{equation}
are dual to $e_1, e_2, e_3$ and
\begin{equation}
de^1=0, \quad de^2=0,
\quad d e^3=d(dz-xdy)=-dx\wedge dy=- e^1 \wedge e^2.
\end{equation}
Now we are going to consider examples of solvmanifolds that
are not nilmanifolds.

3) let $G_1$ be a solvable Lie group of matrices
\begin{equation}
\begin{pmatrix}e^{kz}&0&0&x\\
0&e^{{-}kz}&0&y\\
0&0&1&z\\
0&0&0&1
\end{pmatrix},
\end{equation}

where
$e^k+e^{{-}k}=n \in {\mathbb N}, k \ne 0$.

$G_1$ can be regarded as a
semidirect product
$G_1={\mathbb R} \ltimes {\mathbb R}^2$ where ${\mathbb R}$ acts
on ${\mathbb R}^2$ (with coordinates $x,y$) via
$$
z \to \phi(z)=\begin{pmatrix}e^{kz}&0\\
0 & e^{{-}kz} \end{pmatrix}.
$$
A lattice $\Gamma_1$ in $G_1$ is generated by the following matrices:
$$
\begin{pmatrix}e^{k}&0&0&0\\
0&e^{{-}k}&0&0\\
0&0&1&1\\
0&0&0&1
\end{pmatrix},\quad
\begin{pmatrix}1&0&0&u_1\\
0&1&0&v_1\\
0&0&1&0\\
0&0&0&1
\end{pmatrix}, \quad
\begin{pmatrix}1&0&0&u_2\\
0&1&0&v_2\\
0&0&1&0\\
0&0&0&1
\end{pmatrix},
$$
where
$\begin{vmatrix} u_1&v_1\\ u_2&v_2
\end{vmatrix} \ne 0$.

The corresponding Lie algebra $\mathfrak{g_1}$ has the following basis:
$$
e_1=\begin{pmatrix} k&0&0&0\\
0&{-}k&0&0\\
0&0&0&1\\
0&0&0&0
\end{pmatrix},\quad e_2=
\begin{pmatrix}0&0&0&1\\
0&0&0&0\\
0&0&0&0\\
0&0&0&0
\end{pmatrix}, \quad e_3=
\begin{pmatrix}0&0&0&0\\
0&0&0&1\\
0&0&0&0\\
0&0&0&0
\end{pmatrix},
$$
and the following structure relations:
$$[e_1,e_2]=ke_2, \quad [e_1,e_3]=-ke_3, \quad [e_2,e_3]=0.$$
The left-invariant $1$-forms
\begin{equation}
e^1=dz, \quad e^2=e^{{-}kz}dx, \quad e^3=e^{kz}dy
\end{equation}
are the dual basis to $e_1, e_2, e_3$ and
\begin{equation}
de^1=0, \quad de^2={-}ke^{{-}kz} dz \wedge dx ={-}ke^1 \wedge e^2,
\quad d e^3=k e^1 \wedge e^3.
\end{equation}
As the solvable Lie group $G$ is simply connected the fundamental
group $\pi_1(G/\Gamma)$ is naturally isomorphic to the lattice $\Gamma$:
$\pi_1(G/\Gamma)\cong \Gamma$.

The Lie algebra $\mathfrak{g_1}$ of $G_1$ considered above
is an example of completely solvable Lie algebra.
A Lie algebra $\mathfrak{g}$ is called completely solvable
if $\forall X \in \mathfrak{g}$ operator $ad(X)$
has only real eigenvalues.

Let $G/\Gamma$ be a solvmanifold.
One can identify its de Rham complex $\Lambda^*(G/\Gamma)$
with the subcomplex in  $\Lambda^*(G)$
$$\Lambda^*_{\Gamma inv}(G) \subset \Lambda^*(G)$$
of left-invariant forms on $G$ with respect to the action of the lattice
$\Gamma$.

The subcomplex $\Lambda^*_{\Gamma inv}(G)$ containes in its turn the
subcomplex $\Lambda^*_{\Gamma inv}(G)$ of left-invariant forms
with respect to the action of $G$.

Taking left-invariant vector fields $X_1, {\dots}, X_{p{+}1}$
and a left-invariant $p$-form $\xi\in \Lambda^*_{G inv}(G)$
in the formula (\ref{Cartan_form}) we have:
\begin{equation}
\label{d_leftinv}
d\xi(X_1, {\dots}, X_{p{+}1})=
\sum_{1{\le}i{<}j{\le}p{+}1}({-}1)^{i{+}j}
\xi([X_i,X_j],X_1, {\dots}, \hat X_i,
{\dots}, \hat X_j, {\dots}, X_{p{+}1}).
\end{equation}
The Lie algebra of left-invariant vector fields on $G$
is naturally isomorphic to the tangent Lie algebra $\mathfrak{g}$.
Hence one can identify the space $\Lambda^p_{G inv}(G)$ with
the space $\Lambda^p (\mathfrak{g}^*)$
of skew-symmetric polylinear functions on $\mathfrak{g}$.

The differential $d$ defined by (\ref{d_leftinv}) provides us
with the cochain complex of the Lie algebra $\mathfrak{g}$:
\begin{equation}
\begin{CD}
{\mathbb R} @>{d_0{=}0}>>
\mathfrak{g}^* @>{d}>> \Lambda^2 (\mathfrak{g}^*) @>{d}>>
\Lambda^3 (\mathfrak{g}^*) @>{d}>>
\dots \end{CD}
\end{equation}
The dual of the Lie bracket
$[,]: \Lambda^2 (\mathfrak{g}) \to \mathfrak{g}$ gives a linear mapping
$$\delta: \mathfrak{g}^* \to \Lambda^2 (\mathfrak{g}^*).$$
Consider a basis $e_1,\dots, e_n$ of $\mathfrak{g}$
and its dual basis $e^1,\dots, e^n$. Then we have the following relation:
\begin{equation}
\label{de^k}
de^k=-\delta e^k=-\sum_{i{<}j} C^k_{ij}de^i \wedge de^j,
\end{equation}
where $[e_i,e_j]=\sum C^k_{ij}e_k$.
The differential $d$ is completely determined by (\ref{de^k}) and the
following property:
$$
d(\xi_1 \wedge \xi_2)=d\xi_1 \wedge \xi_2+(-1)^{deg\xi_1} \xi_1 \wedge d\xi_2,
\; \forall \xi_1, \xi_2 \in \Lambda^{*} (\mathfrak{g}^*).
$$
Cohomology of the complex $(\Lambda^*(\mathfrak{g}^*),\delta)$ is called
the cohomology (with trivial coefficients) of the Lie algebra
$\mathfrak{g}$ and is denoted by $H^*(\mathfrak{g})$.

Let us consider the inclusion
$$\psi: \Lambda^*(\mathfrak{g}) \to
\Lambda^*(G/\Gamma).$$

Let $G/\Gamma$ be a compact solvmanifold, where
$G$ is a completely solvable Lie group, then
$\psi: \Lambda^*(\mathfrak{g}) \to
\Lambda^*(G/\Gamma)$ induces the isomorphism
$\psi^*: H^*(\mathfrak{g}) \to H^*(G/\Gamma, {\mathbb R})$
in cohomology (Hattori's theorem \cite{H}, Nomizu's theorem for
nilmanifolds \cite{Nz}).

Let us return to our examples:

1) the cohomology classes $H^*({\mathbb T}^n, {\mathbb R})$
are represented by invariant forms
$$dx^{i_1}\wedge \dots \wedge dx^{i_q}, \quad 1\le i_1 <\dots <i_q\le n,
\; q=1,\dots,n;$$

2) $H^*({\mathcal H}_3/\Gamma_3, {\mathbb R})$ is spanned
by the cohomology classes of the following left-invariant forms:
$$dx, \;dy, \; dy \wedge dz, \; dx\wedge (dz-xdy),\;
dx\wedge dy \wedge dz.$$

3) $H^*(G_1/\Gamma_1, {\mathbb R})$ is spanned
by the cohomology classes of:
$$e^1=dz, \; e^2\wedge e^3= dx \wedge dy, \;
e^1 \wedge e^2\wedge e^3=dx\wedge dy \wedge dz.$$

\section{Deformed differential and Lie algebra cohomology}

From the definition of Lie algebra cohomology
it follows that
$H^1(\mathfrak{g})$ is the dual space to
$\mathfrak{g}/[\mathfrak{g},\mathfrak{g}]$.

1) $b^1(\mathfrak{g})=\dim H^1(\mathfrak{g}) \ge 2$ for
a nilpotent Lie algebra $\mathfrak{g}$ (Dixmier's theorem \cite{D});

2) $b^1(\mathfrak{g}) \ge 1$ for
a solvable Lie algebra $\mathfrak{g}$;

3) $b^1(\mathfrak{g})=0$
for a semi-simple Lie algebra $\mathfrak{g}$.

Consider a Lie algebra $\mathfrak{g}$
with a non-trivial $H^1(\mathfrak{g})$.
Let $\omega \in \mathfrak{g}^*, d \omega =0$. One can define

1) a new deformed differential $d_{\omega}$
in $\Lambda^{*} (\mathfrak{g}^*)$
by the formula
$$d_{\omega}(a)=da+ \omega \wedge a.$$

2) a one-dimensional representation
$$\rho_{\omega}: \mathfrak{g} \to {\mathbb K}, \:
\rho_{\omega}(\xi)= \omega(\xi), \xi \in \mathfrak{g}.$$

Now we recall the definition of Lie algebra cohomology associated
with a representation. Let $\mathfrak{g}$ be a Lie algebra and
$\rho: \mathfrak{g} \to \mathfrak{gl}(V)$ its linear representation.
We denote by $C^q(\mathfrak{g},V)$
the space of $q$-linear alternating mappings of $\mathfrak{g}$ into
$V$. Then one can consider an algebraic complex:
$$
\begin{CD}
V=C^0(\mathfrak{g}, V) @>{d}>>
C^1(\mathfrak{g}, V) @>{d}>> C^2(\mathfrak{g}, V) @>{d}>>
C^3(\mathfrak{g}, V) @>{d}>>
\dots
\end{CD}
$$
where the differential $d$ is defined by:
\begin{equation}
\begin{split}
(d f)(X_1, \dots, X_{q{+}1})=
\sum_{i{=}1}^{q{+}1}(-1)^{i{+}1}
\rho(X_i)(f(X_1, \dots, \hat X_i, \dots, X_{q{+}1}))+\\
+ \sum_{1{\le}i{<}j{\le}q{+}1}(-1)^{i{+}j{-}1}
f([X_i,X_j],X_1, \dots, \hat X_i, \dots, \hat X_j, \dots, X_{q{+}1}).
\end{split}
\end{equation}
The cohomology of the complex $(C^*(\mathfrak{g}, V), d)$ is called
the cohomology of the Lie algebra $\mathfrak{g}$
associated to the representation $\rho: \mathfrak{g} \to \mathfrak{gl}(V)$.

Let $\mathfrak{g}$ be a Lie algebra and $\omega \in \mathfrak{g^*}$
is a closed $1$-form.
Then the complex $(\Lambda^*(\mathfrak{g^*}), d_{\omega})$
coincides with the cochain complex of the Lie algebra
$\mathfrak{g}$ associated with the one-dimensional representation
$\rho_{\omega}: \mathfrak{g} \to {\mathbb K}$,
where $\rho_{\omega}(\xi)= \omega(\xi), \xi \in \mathfrak{g}$.

The proof follows from the formula:
$$(\omega \wedge a) (X_1, \dots, X_{q{+}1})=
\sum_{i{=}1}^{q{+}1}(-1)^{i{+}1}
\omega(X_i)(a(X_1, \dots, \hat X_i, \dots, X_{q{+}1})).
$$

The cohomology $H^{*}_{\omega}(\mathfrak{g})$
coincides with the Lie algebra cohomology with trivial coefficients if
$\omega=0$. If $\omega \ne 0$
the deformed differential $d_{\omega}$
is not compatible with the exterior product $\wedge$
in $\Lambda^*(\mathfrak{g})$
$$d_{\omega}(a \wedge b)=d(a \wedge b)+
\omega \wedge a \wedge b
\ne d_{\omega}(a) \wedge b +
(-1)^{deg a}a \wedge d_{\omega}(b) $$
and the cohomology $H^{*}_{\omega}(\mathfrak{g})$ has no
natural multiplicative structure.

Let $G/\Gamma$ be a compact solvmanifold,
where $G$ is a completely solvable Lie group and
$\tilde \omega$ is a closed $1$-form on $G/\Gamma$.
From the previous sections it follows that
the cohomology $H^*_{\tilde \omega}(G/\Gamma, {\mathbb C})$
is isomorphic to the Lie algebra cohomology
$H^*_{\omega}(\mathfrak{g})$ where
$\omega \in \mathfrak{g}^*$ is the left-invariant $1$-form that represents
the class $[\tilde \omega] \in H^1(G/\Gamma, {\mathbb R})$.

One can define by means of $\omega$
a one-dimenisional representation
$\rho_{\omega}: G \to {\mathbb C}^*$:
$$ \rho_{\omega}(g)=\exp \int_{\gamma(e,g)} \omega,$$
where $\gamma(e,g)$ is a path connecting the identity $e$ with
$g \in G$ (let us recall that $G$ is a symply connected).
As $\omega$ is the left invariant $1$-form then
$$ \int_{\gamma(e,g_1 g_2)} \omega=
\int_{\gamma(e,g_1)} \omega+
\int_{\gamma(g_1,g_1 g_2)} \omega=
\int_{\gamma(e,g_1)} \omega+
\int_{g_1^{{-}1}\gamma(e,g_2)} \omega
$$
holds on and
$\rho_{\omega}(g_1g_2)=\rho_{\omega}(g_1)
\rho_{\omega}(g_2)$.

The representation
$\rho_{\omega}$
induces the representation of corresponding Lie algebra
$\mathfrak{g}$ (we denote it by the same symbol):
$\rho_{\omega}(X)= \omega(X)$.

Let $\mathfrak{g}$ be a $n$-dimensional
real completely solvable Lie algebra (or complex solvable) and
$b^1(\mathfrak{g})=\dim H^1(\mathfrak{g})=k \ge 1$.
Then exists a basis $e^1, \dots, e^n$
in $\mathfrak{g}^*$ such that
\begin{equation}
\begin{split}
\label{solvemodel}
d e^1= \dots = d e^k=0,\\
d e^{k{+}s} = \alpha_{k{+}s} \wedge e^{k{+}s} +
P_{k{+}s}(e^1,\dots, e^{k{+}s{-}1}), \; \;s=1,\dots, n-k.
\end{split}
\end{equation}
where
\begin{equation}
\begin{split}
\alpha_{k{+}s}=\alpha_{s{;}1}e^1+\alpha_{s{;}2}e^2+\dots+
\alpha_{s{;}k}e^k,\\
P_{k{+}s}(e^1,\dots, e^{k{+}s{-}1})=
\sum\limits_{1{\le}i{<}j{\le}k{+}s{-}1} P_{s;i,j}
e^i \wedge e^j.
\end{split}
\end{equation}
It is convinient to define $\alpha_i=0, i=1,\dots,k$.
The set $\{\alpha_1, \dots, \alpha_{n}\}$ of closed $1$-forms is in fact
the set of the weights of completely reducible representation
associated to the adjoint representation $X \to ad(X)$.

For the proof we apply Lie's theorem to
the adjoint representation $ad$ restricted
to the commutant $[\mathfrak{g}, \mathfrak{g}]$:
$$X {\in} \mathfrak{g}\to ad(X):[\mathfrak{g}, \mathfrak{g}]
\to [\mathfrak{g}, \mathfrak{g}].$$
Namely we can choose a basis $e_{k{+}1}, \dots, e_{n}$
in $[\mathfrak{g}, \mathfrak{g}]$ such that
the subspaces $V_i, i=k{+}1, {\dots}, n$ spanned by $e_{i}, \dots, e_{n}$
are invariant with respect to the representation $ad$.
Then we add $e_{1}, \dots, e_{k}$
in order to get a basis of the whole $\mathfrak{g}$.
For the forms of the dual basis $e^{1}, \dots, e^{n}$
in $\mathfrak{g}^*$ we have formulas ~(\ref{solvemodel}).

Let us consider a new canonical basis of $\mathfrak{g}^*$:
\begin{equation}
\begin{split}
\tilde e^1=e^1, \dots , \tilde e^k=e^k, \\
\tilde e^{k{+}s}=t^{2(s{-}1)}e^{k{+}s},
\;\;s=1,\dots,n-k.
\end{split}
\end{equation}
where $t>0$ is a real parameter.

Then for the differential $d_{\omega}$ in the complex
$\Lambda^*(\tilde e^1,\dots, \tilde e^n)$
we have:
$$d_{\omega}=d_0+ \omega \wedge+td_1+t^2d_2+\dots,
\quad d_0 \tilde e^i=\alpha_i \wedge \tilde e^i.$$
In particular
$$
(d_0+ \omega\wedge)(\tilde e^{i_1}\wedge\dots \wedge \tilde e^{i_q})=
(\alpha_{i_1}+\dots+\alpha_{i_q}+ \omega)\wedge
\tilde e^{i_1}\wedge\dots \wedge \tilde e^{i_q}.$$
Now one can define the scalar product in
$\Lambda^q(\tilde e^1,\dots, \tilde e^n)$
declaring the set $\{ e^{i_1}\wedge\dots \wedge e^{i_q} \}$
of basic $q$-forms as an orthonormal basis of
$\Lambda^q(\tilde e^1,\dots, \tilde e^n)$.
Then
\begin{equation}
\begin{split}
d_{\omega}^*d_{\omega}+
d_{\omega}d_{\omega}^*=
R_0+tR_1+t^2R_2+\dots,\\
R_0(\tilde e^{i_1}\wedge\dots \wedge \tilde e^{i_q})=
\|\alpha_{i_1}+\dots+\alpha_{i_q}+ \omega\|^2
\tilde e^{i_1}\wedge\dots \wedge \tilde e^{i_q}.
\end{split}
\end{equation}
As $t \to 0$ the minimal eigenvalue of
$d_{\omega}^*d_{\omega}+
d_{\omega}d_{\omega}^*$ converges to the
minimal eigenvalue of $R_0$. Thus if
$$\alpha_{i_1}+\dots+\alpha_{i_q}+ \omega \ne 0, \quad
1 \le i_1 < i_2 < \dots < i_q \le n $$
then $H^q_{\omega}(\mathfrak{g})=0$.

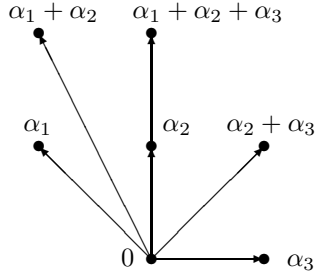
\begin{figure}
\begin{picture}(100,36)(-15,2)
  \put(30,1){\vector(-1,1){15}}
  \put(30,1){\vector(1,0){15}}
  \put(30,1){\vector(0,1){15}}
  \put(30,1){\vector(1,1){15}}
  \put(30,1){\vector(-1,2){15}}
  \put(30,1){\vector(0,1){30}}
  \put(30,1){\circle*{1,5}}
  \put(30,16){\circle*{1,5}}
  \put(45,16){\circle*{1,5}}
  \put(15,16){\circle*{1,5}}
  \put(15,31){\circle*{1,5}}
  \put(30,31){\circle*{1,5}}
  \put(45,1){\circle*{1,5}}
  \put(26,0){$0$}
  \put(31.5,18){$\alpha_2$}
  \put(40,18){$\alpha_2+\alpha_3$}
  \put(13,18){$\alpha_1$}
  \put(11,33){$\alpha_1+\alpha_2$}
  \put(27.5,33){$\alpha_1+\alpha_2+\alpha_3$}
  \put(48,0){$\alpha_3$}
            \end{picture}
\caption{The finite subset $\Omega_{\mathfrak{g}} \subset H^1(\mathfrak{g})$.}
\label{fig1}
\end{figure}

Recall that $\alpha_1=\dots=\alpha_k=0$ and let
us introduce the finite subset
$\Omega_{\mathfrak{g}} \subset H^1(\mathfrak{g})$
such that:
\begin{equation}
\Omega_{\mathfrak{g}}=\left\{ \alpha_{i_1}{+} \dots{+} \alpha_{i_s} |
\;\; 1 \le i_1{<} \dots {<} i_s \le n, \;\; s=1, {\dots}, n \right\}.
\end{equation}
It follows that if
$$-\omega \notin \Omega_{\mathfrak{g}}$$
then the total cohomology $H^*_{\omega}(\mathfrak{g})$ is trivial:
$H^*_{\omega}(\mathfrak{g}) \equiv 0$.

One can easily remark that the subset $\Omega_{\mathfrak{g}}$
is well-defined and does not depend on the ordering of weights $\alpha_i$.

Let $G/\Gamma$ be a compact solvmanifold, where
$G$ is a completely solvable Lie group.
Then the left-invariant closed $1$-forms from
$\Omega_{\mathfrak{g}}$ define a finite subset in
$H^1(G/\Gamma, {\mathbb R})$. We denote this subset by $\Omega_{G/\Gamma}$.
Let
$\omega$ be a closed $1$-form on $G/\Gamma$.
If the cohomology class $$-[\omega] \notin
\Omega_{G/\Gamma}$$ then
the total cohomology $H^*_{\omega}(G/\Gamma, {\mathbb R})$ is trivial:
$H^*_{\omega}(G/\Gamma, {\mathbb R}) \equiv 0$.
The subset $\Omega_{G/\Gamma}$ is well-defined in terms of the
corresponding Lie algebra $\mathfrak{g}$. The corresponding Lie algebra
$\mathfrak{g}$ must to be unimodular, i.e.
the left-invariant $n$-form $e^1\wedge \dots \wedge e^n$
determines non-exact volume form on $G/\Gamma$ and hence
$$\alpha_1+\alpha_2+\dots+\alpha_n=0.$$
If $G/\Gamma$ is a compact nilmanifold then
all the weights $\alpha_i, i=1,\dots, n$ are trivial and
therefore $\Omega_{G/\Gamma}=\{0\}$.
Hence the cohomology $H^*_{\omega}(G/\Gamma, {\mathbb R})$
of a nilmanifold $G/\Gamma$ is trivial
if and only if the form $\omega$ is non-exact.

Let us consider a $3$-dimensional solvmanifold $G_1/\Gamma_1$
defined in the previous section.
We recall that the corresponding Lie algebra $\mathfrak{g}_1$
is defined by its basis $e_1, e_2, e_3$
and the following nontrivial brackets:
$$[e_1,e_2]=k e_2, \quad [e_1,e_3]=-k e_3.$$
For the dual basis of left-invariant $1$-forms
$e^1=dz, e^2=e^{{-}kz}dx, e^3=e^{kz}dy$ we had
$$ de^1=0,\; de^2=-ke^1 \wedge e^2,
\; de^3=k e^1 \wedge e^3.$$
Hence $\alpha_1=0, \alpha_2=-ke^1,
\alpha_3=ke^1$ and $\alpha_2+\alpha_3=0$.

So it is easy to see that
$$\Omega_{G_1/\Gamma_1}=\{\:\pm \: k[e^1]\}$$
and therefore the cohomology
$H^*_{\omega}(G_1/\Gamma_1,{\mathbb R})$ is trivial
if $[\omega] \ne 0, \pm k[e^1]$.

a) $H^*_{k[e^1]}(G_1/\Gamma_1,{\mathbb R})$ is spanned by two classes:
$$e^2=e^{{-}kz}dx, \quad e^1 \wedge e^2=dz \wedge e^{{-}kz}dx.$$

b) $H^*_{{-}k[e^1]}(G_1/\Gamma_1,{\mathbb R})$ is spanned by two classes:
$$e^3=e^{kz}dy, \quad e^1 \wedge e^3=dz \wedge e^{kz}dy.$$

\begin{figure}
\begin{picture}(100,10)(-15,2)
  \put(0,1){\line(1,0){100}}
  \put(30,1){\vector(-1,0){20}}
  \put(30,1){\vector(1,0){20}}
  \put(80,1){\vector(1,0){20}}
  \multiput(9,1)(21,0){3}{\circle*{1,5}}
  \put(4.5,3){$-k[e^1]$}
  \put(29.2,3){$0$}
  \put(49,3){$k[e^1]$}
  \put(71,3){$H^1(G_1/\Gamma_1,{\mathbb R})={\mathbb R}$}
          \end{picture}
\caption{The finite subset $\Omega_{G_1/\Gamma_1}$ }
\label{fig2}
\end{figure}
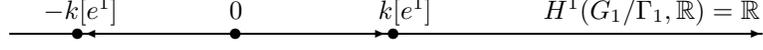

Hence we have the following Betti numbers
$b^p_{\omega}{=}
\dim H^p_{\omega}(G_1/\Gamma_1, {\mathbb R})$
of the solvmanifold $G_1/\Gamma_1$:
\begin{equation}
\begin{split}
1) \quad b_{\pm ke^1}^0=0, \; b^1_{\pm ke^1}=
b^2_{\pm k e^1}=1,\; b^3_{\pm k e^1}=0;\\
2) \quad b^0_0 = b^1_0=
b^2_0 = b^3_0= 1.
\end{split}
\end{equation}

It was proved by G.~Mostow in ~\cite{Mos} that
any compact solvmanifold $G/\Gamma$ is a bundle with toroid as base space
and nilmanifold as fibre, in particular a solvmanifold $G/\Gamma$
is fibered over the circle $\pi: G/\Gamma \to S^1$. Hence the $1$-form
$\pi^*(d\varphi)$ on $G/\Gamma$ has no critical points:
$m_p(\pi^*(d\varphi))=0,\; \forall p$.
It follows from A.~Pajitnov's theorem \cite{Pa}
that for $\lambda$ sufficiently
large we have
$H^p_{\lambda \pi^*(d\varphi)}(G/\Gamma,{\mathbb R})=0, \; \forall p$.

Now we are going introduce an example of solvmanifold $G/\Gamma$
with non completely solvable Lie group $G$.
Let $G_2$ be a solvable Lie group of matrices
\begin{equation}
\begin{pmatrix}\cos{2\pi z}&\sin{2\pi z}&0&x\\
{-}\sin{2\pi z}&\cos{2\pi z}&0&y\\
0&0&1&z\\
0&0&0&1
\end{pmatrix}.
\end{equation}

A lattice $\Gamma_2$ in $G_2$ is generated by the following matrices:
$$
\begin{pmatrix}\cos{\frac{2\pi n}{p}}&\sin{\frac{2\pi n}{p}}&0&0\\
{-}\sin{\frac{2\pi n}{p}}&\cos{\frac{2\pi n}{p}}&0&0\\
0&0&1&\frac{n}{p}\\
0&0&0&1
\end{pmatrix}, \quad
\begin{pmatrix}1&0&0&u_1\\
0&1&0&v_1\\
0&0&1&0\\
0&0&0&1
\end{pmatrix}, \quad
\begin{pmatrix}1&0&0&u_2\\
0&1&0&v_2\\
0&0&1&0\\
0&0&0&1
\end{pmatrix},
$$
where $n$ is an integer, $p=2,3,4,6$ and
$\begin{vmatrix} u_1&v_1\\ u_2&v_2
\end{vmatrix} \ne 0$, or another type:
$\tilde \Gamma_2$ is generated by the following
matrices:
$$
\begin{pmatrix}1&0&0&1\\
0&1&0&0\\
0&0&1&0\\
0&0&0&1
\end{pmatrix}, \quad
\begin{pmatrix}1&0&0&0\\
0&1&0&1\\
0&0&1&0\\
0&0&0&1
\end{pmatrix},
\begin{pmatrix}1&0&0&u\\
0&1&0&v\\
0&0&1&n\\
0&0&0&1
\end{pmatrix}, \quad
$$
where $n$ is an integer.
The corresponding Lie algebra $\mathfrak{g_2}$ has the following basis:
$$
e_1=\begin{pmatrix} 0&2\pi&0&0\\
{-}2\pi&0&0&0\\
0&0&0&1\\
0&0&0&0
\end{pmatrix},\quad e_2=
\begin{pmatrix}0&0&0&1\\
0&0&0&0\\
0&0&0&0\\
0&0&0&0
\end{pmatrix}, \quad e_3=
\begin{pmatrix}0&0&0&0\\
0&0&0&1\\
0&0&0&0\\
0&0&0&0
\end{pmatrix},
$$
and the following structure relations:
$$[e_1,e_2]=-2\pi e_3, \quad [e_1,e_3]=2\pi e_2, \quad [e_2,e_3]=0.$$
As the eigenvalues of $ad(e_1)$ are equal to $0,\pm 2\pi i$
the Lie group $G_2$ is not completely solvable.

The left-invariant $1$-forms
\begin{equation}
e^1=dz, \quad e^2=\cos{2\pi z}dx-\sin{2\pi z}dy,
\quad e^3=\sin{2\pi z}dx+\cos{2\pi z}dy
\end{equation}
are the dual basis to $e_1, e_2, e_3$ and
\begin{equation}
de^1=0, \quad de^2={-}2 \pi e^1 \wedge e^3,
\quad d e^3=2 \pi e^1 \wedge e^2.
\end{equation}
The cohomology $H^*(\mathfrak{g}_2)$ is spanned
by the cohomology classes of:
$$e^1, \; e^2\wedge e^3, \; e^1 \wedge e^2\wedge e^3.$$
But $$\dim{H^1(\mathfrak{g}_2)}=1 \ne
\dim{H^1(G_2/\Gamma_2, {\mathbb R})}=3.$$
This example shows that, generally speaking, Hattori's theorem
does not hold for non completely solvable Lie groups $G$,
but the inclusion of left-invariant differential forms
$\psi: \Lambda^*(\mathfrak{g}^*) \to \Lambda^*(G/\Gamma)$
always induces the injection
$\psi^*$ in cohomology.

\end{document}